\documentclass[seceqn,secthm]{elsart}

\usepackage{amsfonts}
\usepackage[dvips]{graphics,color}
\newcommand{\Z}{\mathbb{Z}}
\newcommand{\Q}{\mathbb{Q}}

\newcommand{\A}{\mathbb{A}}

\input amssym.def
\newsymbol\square 1003
\newsymbol\lozenge 1006
\newsymbol\surjarrow 1310
\newsymbol\injarrow 131A
\newsymbol\ltimes 226E
\newsymbol\rtimes 226F
\newsymbol\void 203F
\newsymbol\leftrightsquigarrow 1321
\def\END {$\square$}

\def\proclaim #1. #2\par #3\par {\medbreak
\noindent{\bf#1.\enspace}{\sl#2}\par\medbreak
\noindent{\bf Proof.} #3 \par}
\def\proclaimb #1. #2\par {\medbreak
\noindent{\bf#1.\enspace}{\sl#2}\par\medbreak}
 0
 2
 0
\begin{document}
\begin{frontmatter}
\title{Homology of holomorphs of free groups}
\author{Craig A. Jensen$^1$}
\thanks{Partially supported by Louisiana Board of Regents Research Competitiveness Subprogram Contract LEQSF-RD-A-39.}

\smallskip

\address{Department of Mathematics, University of New Orleans\\
New Orleans, LA 70148, USA\\
{\tt jensen@math.uno.edu}}

\begin{abstract}
The holomorph of a free group $F_n$ is the semidirect product
$F_n \rtimes Aut(F_n)$.  Using the methods of Hatcher and Vogtmann
in \cite{[H-V]} and \cite{[V]}, we derive stability results and
calculate the mod-$p$ homology of 
these holomorphs for odd primes $p$ in dimensions 1 and 2, and their rational
homology in dimensions 1 through 5.  Calculations of
the twisted (where $Aut(F_n)$ acts by first projecting to $Gl_n(\Z)$
and then including in $Gl_n(\Q)$)
homology $H_*(Aut(F_n); \Q^n)$ follow in corresponding dimensions.
\end{abstract}

\end{frontmatter}

\section{Introduction}

Let $F_n$ denote the free group on $n$ letters and let
$Aut(F_n)$ and $Out(F_n)$ denote the automorphism group
and outer automorphism group, respectively, of $F_n$.
Define the {\em holomorph} of $F_n$ to be the semidirect
product $F_n \rtimes Aut(F_n)$.

When examining the action of $Aut(F_n)$ on auter space 
(defined in \cite{[H]} and \cite{[H-V]} by Hatcher and Vogtmann) 
or the action of the symmetric automorphism group $\Sigma Aut(F_n)$
on a corresponding space (see Collins in \cite{[C1]} and 
McCullough and Miller in \cite{[M-M]}), holomorphs arise 
naturally by looking at point stabilizers.  Hence (see \cite{[J1]}
and \cite{[JD]}, where results like these are used) 
calculating the homology of lower rank holomorphs
of $F_n$ is very useful when calculating the homology of 
$Aut(F_n)$ (or related groups, like $\Sigma Aut(F_n)$) in higher 
ranks.  Holomorphs of free groups have also been studied 
recently by Thomas and Velickovic in \cite{[T-V]}.

Because the Hochschild-Serre spectral sequence of the group extension 
$1 \to F_n \to F_n \rtimes Aut(F_n) \to Aut(F_n) \to 1$ has
$E^2$-page 
$$E^2_{r,s} = H_r(Aut(F_n); H_s(F_n; M)) \Rightarrow 
H_{r+s}(F_n \rtimes Aut(F_n); M),$$ we see that
calculating the homology of the holomorph also amounts to 
calculating the twisted homology of $Aut(F_n)$ with coefficients 
in $M^n$.
See Dwyer \cite{[D]} or Borel \cite{[Bo]} for results about
$H_*(GL_n(\Z); (\Z/p)^n)$,
$H_*(GL_n(\Z); \Q^n)$ which motivated this paper.  See also
Allison, Ash, and Conrad \cite{[A]} and Charney \cite{[Ch]} for other 
results in the field.

Let $F_n$ have free basis $\{x_1, \ldots, x_n\}$.  Define a
preferred inclusion $\iota: F_n \to F_{n+1}$ on the
generators by setting $\iota(x_i)=x_i$ for all $1 \leq i \leq n$.
Note that $\iota$ includes $F_n$ as a free factor of $F_{n+1}$.
There is an induced preferred inclusion $\iota: Aut(F_n) \to Aut(F_{n+1})$
defined by $\iota(\phi)(x_i)=\phi(x_i)$ if $1 \leq i \leq n$
and $\iota(\phi)(x_{n+1})=x_{n+1}$.  Finally, there induces
a preferred inclusion 
$\iota: F_n \rtimes Aut(F_n) \to F_{n+1} \rtimes Aut(F_{n+1})$
defined by $\iota(x,\phi) = (\iota(x),\iota(\phi))$.

In this paper, we prove the following two results:

\begin{thm}[Homology stability]  \label{stab}

\noindent \begin{enumerate}
\item The map $H_i(F_n \rtimes Aut(F_n); \Q) \to
H_i(F_{n+1} \rtimes Aut(F_{n+1}); \Q)$ induced by preferred inclusion is
an isomorphism for $n > 3i/2.$
\item The map $H_i(F_n \rtimes Aut(F_n); \Z) \to
H_i(F_{n+1} \rtimes Aut(F_{n+1}); \Z)$ induced by preferred inclusion is
an isomorphism for $n \geq 4i+2.$
\end{enumerate} \end{thm}

\begin{thm}[Low dimensional homology] \label{lowdim}

\noindent \begin{enumerate}
\item If $p$ is an odd prime, $1 \leq i \leq 2$ and $n$ 
is any positive integer, then
$$H_i(F_n \rtimes Aut(F_n); \Z/p) = 0.$$
\item If $1 \leq i \leq 5$ and $n$ is any 
positive integer, then
$$H_i(F_n \rtimes Aut(F_n); \Q) = 0,$$
except when $i=4, n=3$ or $i=4,n=4$ in which case
$$H_i(F_n \rtimes Aut(F_n); \Q) = \Q.$$
\item If $0 \leq i \leq 4$ and $n$ is any 
positive integer, then the twisted homology group
$$H_i(Aut(F_n) ; \Q^n) = 0,$$
except when $i=3, n=3$ where
$$H_3(Aut(F_3) ; \Q^3) = \Q.$$
\end{enumerate}
\end{thm}

Because of the large size of the spaces involved, Maple programs
were used to establish Theorem \ref{lowdim} (2).

In \cite{[H-V]} Hatcher and Vogtmann prove a ``Degree Theorem''
which they use to derive linear stability
ranges for the integral cohomology of
$Aut(F_n)$ and to show in \cite{[V]} that 
for $1 \leq t \leq 6$ and all $n \geq 1$,
$$\matrix{
\hfill H_t(Aut(F_n); \Q) &=& &\left\{\matrix{
\Q \hfill &\hbox{if } t = n = 4 \hfill \cr
0 \hfill &\hbox{otherwise} \hfill \cr} \right. \hfill \cr
}$$
To prove Theorems \ref{stab} and \ref{lowdim}, we borrow the
methods of Hatcher and Vogtmann, and details will be omitted when
we closely follow their work.  

The next section will review auter space and introduce spaces
with natural actions by 
holomorphs of free groups, while the third section 
argues that the degree theorem holds in the new context of holomorphs of 
free groups.  The fourth section is devoted to proving Theorem
\ref{stab},
and the fifth section contains a proof of Theorem \ref{lowdim}.

The author would like to thank
Prof. Vogtmann (his thesis advisor while at Cornell) and Prof. Hatcher
for their help through the years.  The author would like to thank the
referee for many helpful comments which improved the exposition of
this paper.

\section{$Aut(F_n)$, auter space, and holomorphs}

Let $(R_n,*)$ be the $n$-leafed rose, $n$ circles wedged together
at the basepoint $*$.
Recall (see  \cite{[C-V]}, \cite{[H-V]}, and \cite{[S-V1]})
that the spine $X_n$ of ``auter space'' is 
the realization of a poset of
pointed graphs $(\Gamma,*)$ equipped with {\em markings} 
(homotopy equivalences) from the $n$-leafed rose.  
The poset structure derives from forest 
collapes in $\Gamma$ and the $Aut(F_n)$-action twists the 
marking.  Let $Q_n$ be the quotient of $X_n$ by
$Aut(F_n)$.  Define $\A_n$ to be auter space (as opposed to its spine)
where the edges in a marked graph have lengths which must sum to 1.
Define $X_{n,k}$, $Q_{n,k}$, and $\A_{n,k}$ to be the parts of
$X_{n}$, $Q_{n}$, and $\A_{n}$, respectively, corresponding to 
graphs of degree less than or equal to $k$ (see \cite{[H-V]}.)

Let $\Theta_m$ be the graph with two vertices
and $m+1$ edges, each of which goes from one vertex to the other.
Let $\mathcal{Q} \cong \Z/3$ be the subgroup of $Aut(F_{n+2})$
given by the $\Z/3$-action of rotating the edges of the graph
$\Theta_2$ in $R_n \vee \Theta_2$ (where the wedge joins the vertex of
$R_n$ with one of the two vertices of $\Theta_2$ and the
resulting vertex is the basepoint of the graph.)  In 
Definition 6.3 of \cite{[J1]}, 
the space $\tilde X_n$ was defined to be the fixed point
subcomplex $X_{n+2}^\mathcal{Q}$ and it was observed that
the normalizer
$$N_{Aut(F_{n+2})}(\mathcal{Q}) \cong \Sigma_3
\times (F_n \rtimes Aut(F_n))$$
acts properly 
on $X_{n+2}$.
In addition, Definition 6.3 of \cite{[J1]} defines
$\tilde Q_n$ to be the quotient of
$\tilde X_n$ by $N_{Aut(F_{n+2})}(\mathcal{Q})$.
Proposition
6.8 of \cite{[J1]} describes the action 
of $N_{Aut(F_{n+2})}(\mathcal{Q})$ on $\tilde X_n$ in detail.
It is shown that the quotient of $\tilde X_n$ by
$F_n \rtimes Aut(F_n)$ is also $\tilde Q_n$.

As shown in Proposition 6.4 of \cite{[J1]}, the fixed point
space $\tilde X_{n}$ can equivalently be characterized
as the realization of the poset of equivalence classes of
pairs $(\alpha,f)$, where
$\alpha : R_{n} \to \Gamma_{n}$ is
a pointed marked graph whose underlying graph $\Gamma_{n}$
has a special (possibly valence 2)
vertex which is designated as $\circ$,
$\circ$ may equal the basepoint $*$ of $\Gamma_{n}$,
and $f : I \to \Gamma_{n}$ is a homotopy
class (rel endpoints) of maps from
$*$ to $\circ$ in $\Gamma_{n}$.  Basically,
the extra vertex $\circ$ encodes that $\Theta_2$ is
attached at the indicated point when translating back to
marked graphs of genus $n+2$.
Marked graphs
$(\alpha^1,f^1)$ and $(\alpha^2,f^2)$
are {\em equivalent} 
if there
is a homeomorphism $h$ from 
$\Gamma_{n}^1$ to $\Gamma_{n}^2$
which sends $*$ to $*$, $\circ$ to $\circ$,
such that
$$(h \alpha^1)_\# = (\alpha^2)_\# : \pi_1(R_n,*) \to 
\pi_1(\Gamma_n^2)$$
and such that the paths
$$ hf_1, f_2 : I \to \Gamma_n^2$$
are homotopic rel endpoints.
Moreover, from Remark 6.1 of \cite{[J1]}, the quotient
space $\tilde Q_{n}$ 
can be characterized
as the realization of the poset of equivalence classes of
pointed graphs $\Gamma_{n}$
which have possibly valence 2 special vertex $\circ$
(which we think of as signifying that a $\Theta_2$ 
should be attached to the graph at that point)
which may equal the basepoint $*$.

As in Definition 6.9 of \cite{[J1]}, define $\tilde \A_n$ 
to be the analog of $\tilde X_n$ where
the edges in graphs have lengths which must sum to 1.  That is,
in addition to having markings from $\alpha$ and $f$,
graphs in $\tilde \A_n$ have metrics where each edge has a length and the metric is the induced path metric.  We further normalize these metrics
by insisting that the sums of the lengths of the edges of a graph in
$\tilde \A_n$ must sum to 1.

Corresponding to the
 preferred inclusion $Aut(F_n) \to Aut(F_{n+1})$,
there is an equivariant
map of spaces $X_n \to X_{n+1}$ obtained by sending the marked graph
$\alpha: R_n \to \Gamma$ to $(\alpha,id): R_n \vee R_1 \to \Gamma \vee R_1.$
Let $G_n=F_n \rtimes Aut(F_n).$
Similarly, the preferred inclusion 
$G_n \to G_{n+1}$ corresponds to a $G_n$-equivariant preferred inclusion
$\iota: \tilde X_n \to \tilde X_{n+1}$ given by sending the pair
$(\alpha,f)$ to $((\alpha,id),f).$  

\section{The modified degree theorem} \label{c9}

Our goal is to show that Hatcher and Vogtmann's
Degree Theorem in \cite{[H-V]} also applies to holomorphs of
free groups and the space $\tilde X_n$.  This section is
written to convince those already familiar with 
\cite{[H-V]} that their proof carries over into this
new context, and is not meant to be read independently.

As in the case of $\A_n$, graphs in $\tilde \A_n$ come
equipped with a ``height function'' measuring distance
to the basepoint.  Given a particular point $v$ in 
such a graph,
we can take a small neighborhood of it and obtain a star graph
consisting of $v$ and the germs of all edges attached to $v$
in the graph.  Some of these germs are {\em ascending} because points
on them have have height greater than that of $v$ (equivalently,
points on that germ are farther away from the basepoint in the graph.)
Other germs are {\em descending} because points on them are closer 
to the basepoint.

Define the {\em degree} of a marked graph
$$\matrix{\hfill (\alpha,f): R_{n} \coprod I \to \Gamma_{n}\hfill \cr}$$
or, equivalently, the degree of the
underlying graph $\Gamma_{n}$ of the marked graph,
to be
$$deg(\alpha,f) = deg(\Gamma_{n}) = \sum_{v \not = *} ||v|| - 2,$$
where the sum is over all vertices of $\Gamma_{n}$
except the basepoint, and where the {\em augmented valence} 
$||v||$ of a vertex
$v \not = \circ$
is the number of oriented edges starting at $v$ (i.e., the valence 
$|v|$ of $v$.).
The {\em augmented valence} of the vertex $\circ$ is defined to be
one plus the number of oriented edges starting at $\circ$,
or $|\circ|+1$.
This definition is similar to the one given in \cite{[H-V]},
but the intuition (which we do not attempt to
make precise) is that
we treat the vertex $\circ$ as if it signifies that
an ascending germ of an extra edge
is attached to that vertex.
Throughout this section,
our modified definitions of degree, split degree,
canonical splitting, etc., will be motivated by this notion of
thinking that the vertex $\circ$ denotes the germ of another edge
entering that vertex.

Let $\tilde X_{n,k}$, $\tilde Q_{n,k}$, and $\tilde \A_{n,k}$
be the subspaces of
$\tilde X_{n}$, $\tilde Q_{n}$, and $\tilde \A_{n}$,
respectively,
where only marked graphs of degree at most $k$ are considered.
Define $D_k$ to be a $k$-dimensional disk.
We want to prove the following
analog of the degree theorem, Theorem 3.1 of \cite{[H-V]}:

\begin{thm} \label{t13}  A piecewise linear map $f_0:D_k \to \tilde
\A_{n}$ is homotopic to a map $f_1:D_k \to \tilde \A_{n,k}$
by a homotopy $f_t$ during which degree decreases monotonically,
i.e., if $t_1 < t_2$ then $deg(f_{t_1}(s)) \geq deg(f_{t_2}(s))$
for all $s \in D_k$.
\end{thm}

As in \cite{[H-V]}, the immediate corollary is

\begin{cor} \label{t14} The pair $(\tilde \A_{n}, \tilde \A_{n,k})$
is $k$-connected.
\end{cor}

In \cite{[J2]}, we show that $\tilde \A_{n}$ is contractible, so
that Corollary \ref{t14} implies that
$\tilde \A_{n,k}$ is $(k-1)$-connected.

Hatcher and Vogtmann prove the degree theorem by
using various homotopies to deform the
underlying graphs of marked graphs
$$\alpha: R_n \to \Gamma_n.$$
We will use these same homotopies
to deform the marked graphs
$$\matrix{\hfill (\alpha,f): R_{n} \coprod I \to \Gamma_{n}\hfill \cr}$$
that appear in the context of $\tilde \A_n$.
In a remark
following ``Stage 1: Simplifying the critical point'' in
\cite{[H-V]}, Hatcher and Vogtmann
mention that it is obvious
where their homotopies of the
underlying graph $\Gamma_n$ send the basepoint $*$
and the marking $\alpha$.
Our task is to decide where these homotopies send the
extra point $\circ$ on the graph.  It will then be
clear where the homotopies send the path $f$
from $*$ to $\circ$ in the graph $\Gamma_n$.

A point in the graph is a
{\em critical point} if in a small neighborhood of that
point there is more than one
descending germ.

A few notes are in order about specific parts of the
paper by Hatcher and Vogtmann and how these should be
modified:
\begin{enumerate}
\item Canonical splittings.  A procedure called ``canonical
splitting'' is defined in \cite{[H-V]} to decrease
the degrees of graphs in a canonical way.  
A canonical splitting should move
$\circ$ down to the next critical point or the basepoint, provided
$\circ$ is not already a critical point.  
See Figure \ref{holomorph1} for examples.

\bigskip

\centerline{\input{holomorph1.pic}}
\begin{figure}[here]
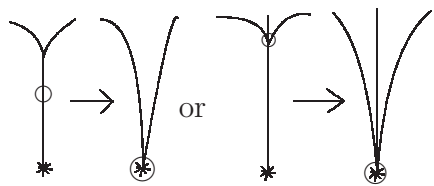

\caption{\label{holomorph1} Canonical splittings}
\end{figure}

\item Sliding $\epsilon$-cones.  We can also perturb the
graph slightly and slide $\circ$ downward off of a 
critical point.

\item Codimension.
As in \cite{[H-V]}, the {\em codimension} of a point on the
graph is one less than the number of downward directions
from that point.  The {\em codimension} of a graph is the
sum of the codimensions of its critical points.

\item Lemma 4.1 needs no modification.

\item Lemma 4.2.  During the homotopies used in Lemma 4.2, only
the lengths of the edges of the underlying graph $\Gamma_{n}$
are perturbed.  The combinatorial
structure of the graph is not changed at
all.  Hence it is clear where $\circ$ and the path $f$ from $*$ to
$\circ$ are sent during these homotopies.  

\item Complexity.
As before, we think of $\circ$ as 
having an ascending germ attached there.  Hatcher and Vogtmann
defined a connecting path as a downward path from one
critical point to another.  We modify this by saying that
the extra germ attached at $\circ$ counts as the beginning of
a downward path that came from a critical point lying above.  
With this convention, the complexity $c_s$ (respectively, $e_s$) 
is defined as the number of connecting paths
(respectively, without critical points in their interiors)
in the graph.  See Figure \ref{holomorph2}.

\centerline{\small \input{holomorph2.pic}}
\begin{figure}[here]
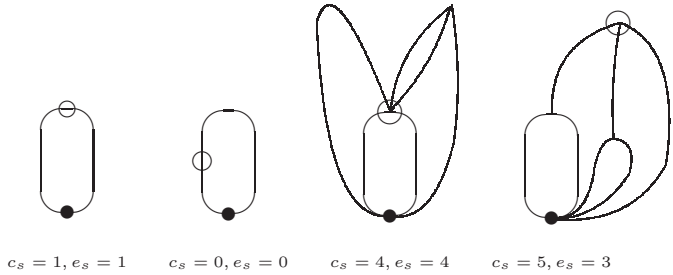

\caption{\label{holomorph2} Complexity examples with $\circ$}
\end{figure}

Using this definition of complexity in place of
that given in \cite{[H-V]}, the section
``canonical splitting and extension to a neighborhood'' 
of \cite{[H-V]} remains valid as written.
For the section
``reducing complexity by sliding in the $\epsilon$-cones'',
consider $\circ$ as giving an attaching point
$\alpha_j$ of a branch $\beta_j$.  Now argue as
directed in \cite{[H-V]}.
\end{enumerate}

Using the above guidelines, the proof of
Theorem \ref{t13} follows from the
work of Hatcher and Vogtmann in \cite{[H-V]}.

\section{Stability results}

The following lemma is a restatement of Lemma 5.2 of \cite{[H-V]}.

\begin{lem} \label{lemma52}
Let $\Gamma$ be the underlying graph of a marked 
graph in $\tilde X_{n,k}$.

\begin{enumerate}
\item If $k < n/2$, then $\Gamma$ has an $R_1$ wedge summand.
\item If $k < 2n/3$, then $\Gamma$ has an $R_1$ or $\Theta_2$
wedge summand.
\item If $k < n-1$, then $\Gamma-\{*\}$ is disconnected.
\end{enumerate}
\end{lem}

\begin{pf*}{Proof.}
Without loss of generality, assume that $\Gamma$ has degree $k$,
that all vertices not equal to $*$ or $\circ$ are trivalent, and
that $\circ$ is bivalent if it is not equal to $*$.  
Define $\Lambda$ to be the full
subgraph of $\Gamma$ spanned by all non-basepoint vertices.
Let $E$ and $V=k$ be the number of edges and vertices,
respectively, of $\Lambda$. 

If $\Gamma$ has no $R_1$ wedge summand at $*$, then the valence $|*|$ is 
$2n-k$ if $*=\circ$ (respectively $2n-k-1$ if $*\not=\circ$) and so
$\chi(\Gamma)=1-n$ is $1-2n+2k-E$ (resp. $-2n+2k-E$); therefore,
$n$ is $2k-E$ (resp. $2k-E-1$) and (1) follows.  If $k<2n/3$, then we
can assume $E$ is less than $k/2$ (resp. $k/2-1$), resulting in an
isolated vertex $v$ of $\Lambda$ not equal to $\circ$ and
establishing (2).  If $k < n-1$, similar arguments show that
$\chi(\Lambda)$ is greater than 1 (resp. 2) so that $\Lambda$ is
disconnected.
\END \end{pf*}

Note that Theorem \ref{stab} (1) follows from 
Lemma \ref{lemma52} in the same way that
Proposition 5.5 of \cite{[H-V]} follows from
Lemma 5.2 of \cite{[H-V]}.  The proof of integral
homology stability contained in \cite{[H-V]} does not
appear, however, to apply to the case of holomorphs.  We try
a different approach here.

\begin{pf*}{Proof of Theorem \ref{stab} (2).} 
Recall that $G_n = F_n \rtimes Aut(F_n)$.
Fix a positive integer $i$ and assume $n \geq 4i+2$. 
Let $C \to \Z$
and $C' \to \Z$ be the augmented cellular chain complexes of
$\tilde X_{n,i+1}$ and $\tilde X_{n+1,i+1}$, respectively.  There
is a chain map $\iota: C \to C'$ induced by the 
preferred inclusion 
$\iota: \tilde X_{n,i+1} \to \tilde X_{n+1,i+1}.$  The 
augmented complexes $C \to Z$ and $C' \to Z$ are exact in dimensions 
less than or equal to $i$ by Corollary \ref{t14}.
Let $F \to \Z$ be
a free resolution of $\Z$ as a $\Z G_{n+1}$-module.  Since
$G_n$ is contained in $G_{n+1}$, $F \to \Z$ is also a free
resolution of $Z$ as a $\Z G_n$-module.  Let $\phi$ be the 
composition
$$F \otimes_{\Z G_n} C \buildrel 1 \otimes \iota \over \to
F \otimes_{\Z G_n} C' \buildrel p \over \to
F \otimes_{\Z G_{n+1}} C'.$$
The morphism $\phi$ of double complexes 
induces morphisms of the spectral sequences corresponding to
the double complexes.  
Taking vertical filtrations of $F \otimes_{\Z G_n} C$
and $F \otimes_{\Z G_{n+1}} C'$ (see \cite{[B]} page 173), we obtain spectral 
sequences where the $E^2_{r,s}$ pages converge to $H_{r+s}(G_n; \Z)$ and
$H_{r+s}(G_{n+1}; \Z)$, respectively, for $r+s \leq i$.  Taking horizontal
filtrations, we have spectral sequences with
$$E_{r,s}^1= \prod_{\sigma \in \tilde Q^r_{n,i+1}} 
H_s(\hbox{stab}_{G_n}(\sigma); \Z) \Rightarrow
H_{r+s}(G_n; \Z), \hbox{ for } r+s \leq i,$$
and
$$\bar E_{r,s}^1= \prod_{\sigma \in \tilde Q^r_{n+1,i+1}} 
H_s(\hbox{stab}_{G_{n+1}}(\sigma); \Z) \Rightarrow
H_{r+s}(G_{n+1}; \Z), \hbox{ for } r+s \leq i,$$
where $\tilde Q^r_{n,i+1}$ and $\tilde Q^r_{n+1,i+1}$ are the
$r$-cells in $\tilde Q_{n,i+1}$ and $\tilde Q_{n+1,i+1}$,
respectively.
Since $n \geq 4i+2$, Lemma \ref{lemma52} (1) yields that 
every marked graph in $\tilde Q_{n,i+1}$ (respectively $\tilde Q_{n+1,i+1}$)
has at least $2i$ (respectively $2i+1$) loops at the basepoint.
This yields a homeomorphism $f: \tilde Q_{n,i+1} \to \tilde Q_{n+1,i+1}$
defined on graphs by adding a loop at the basepoint.
Moreover, if $\sigma \in \tilde Q^r_{n,i+1}$, we can represent
$\sigma$ as a pair $(\Gamma,F)$, where $F$ is some chain of forest
collapses  in $\Gamma$.  
Write $\Gamma=\Gamma_0 \vee R_j$, where
$\Gamma_0$ has no loops at the basepoint and $j \geq 2i$.  Let $G$
be the group of graph isomorphisms of $\Gamma_0$ respecting the
chain $F$ of forests (see \cite{[S-V1]}.) Hence
$f(\sigma)=(\Gamma_0 \vee R_{j+1},F)$,
$\hbox{stab}_{G_n}(\sigma) = G \times \Sigma_j$, and
$\hbox{stab}_{G_{n+1}}(\sigma) = G \times \Sigma_{j+1}$.  From
the stability result for the homology of the symmetric groups in
Corollary 6.7 of \cite{[N]} and the K\"{u}nneth Formula, 
we see that the induced map 
$$H_s(\hbox{stab}_{G_{n}}(\sigma); \Z) \to
H_s(\hbox{stab}_{G_{n+1}}(\sigma); \Z)$$
is an isomorphism for $s \leq i$.  Thus
$\phi_*: E_{r,s}^1 \to \bar E_{r,s}^1$ is an isomorphism for
$r+s \leq i$ and (cf. Proposition 2.6 from Chapter VII
of \cite{[B]}) $\phi_*: H_k(G_n; \Z) \to H_k(G_{n+1}; \Z)$
is an isomorphism for $k \leq i$. 
\END \end{pf*}

\section{Low dimensional homology groups} \label{c10}

As in Hatcher-Vogtmann \cite{[V]}, the 
degree theorem, Theorem \ref{t13},
can be used to prove Theorem \ref{lowdim} and calculate the homology 
$F_n \rtimes Aut(F_n)$ in low dimensions. 

\begin{lem} \label{lemrat} If $n,i$ are positive
integers, then 
$$H_i(F_n \rtimes Aut(F_n); \Q) \cong H_i(\tilde Q_{n,i+1}; \Q)$$
and 
$$H_{i+1}(\tilde Q_{n,i+1}, \tilde Q_{n,i}; \Q) \to
H_i(\tilde Q_{n,i}; \Q)
\to H_i(\tilde Q_{n,i+1}; \Q) \to 0$$ 
is exact.
\end{lem}

\begin{pf*}{Proof.} 
That $H_i(G_n; \Q) \cong H_i(\tilde Q_{n,i+1}; \Q)$
follows from considering the equivariant homology spectral sequence for
$G_n$ acting on $\tilde Q_{n,i+1}$ (cf. \cite{[B]}), noting
that it is concentrated in only one row because we have $\Q$ coefficients
and stabilizers are finite, and since $\tilde Q_{n,i+1}$ is $i$-connected.
For similar reasons, we have that both 
$H_{i+1}(\tilde Q_{n}, \tilde Q_{n,i+1}; \Q)$ and
$H_{i}(\tilde Q_{n}, \tilde Q_{n,i}; \Q)$ are zero so that the long
exact sequence of the triple 
$(\tilde Q_{n}, \tilde Q_{n,i+1}, \tilde Q_{n,i})$ gives that
$H_{i}(\tilde Q_{n,i+1}, \tilde Q_{n,i}; \Q)=0$.  The lemma follows
by considering the long exact sequence of the pair 
$(\tilde Q_{n,i+1},\tilde Q_{n,i}).$
\END \end{pf*}

A corresponding result for $\Z/p$ coefficients is:

\begin{lem} \label{lemmod} Let $p$ be an odd prime and $n \geq 1$.
Then
\begin{enumerate}
\item $H_1(F_n \rtimes Aut(F_n); \Z/p) = 0.$
\item $H_2(F_n \rtimes Aut(F_n); \Z/p) \cong H_2(\tilde Q_{n,3}; \Z/p).$
\item 
$H_{3}(\tilde Q_{n,3}, \tilde Q_{n,2}; \Z/p) \to
H_2(\tilde Q_{n,2}; \Z/p)
\to H_2(\tilde Q_{n,3}; \Z/p) \to 0$
is exact.
\end{enumerate}
\end{lem}

\begin{pf*}{Proof.}
From the five term exact sequence of the group extension corresponding to
the semidirect product $G_n$ (cf. page 171 of \cite{[B]}),
$H_1(F_{n}; \Z)_{Aut(F_{n})} \to H_1(G_n; \Z) \to H_1(Aut(F_{n}); \Z)$
is exact. By \cite{[H]}, $H_1(Aut(F_n); \Z) = \Z/2$.  
Let $\xi_n \in Aut(F_n)$ 
be the automorphism (cf. \cite{[G-J]}) 
sending each generator to its inverse.  Then
$\xi_n$ sends 
$(x_1, \ldots, x_n) \in \Z \oplus \ldots \oplus \Z \cong H_1(Aut(F_n); \Z)$
to $(-x_1, \ldots, -x_n)$.  Hence
$H_1(F_{n}; \Z)_{Aut(F_{n})}$ is all 2-torsion.  Thus $H_1(G_n; \Z)$
is 2-torsion and $H_1(G_n; \Z/p)=0.$  

To show
$H_2(G_n; \Z/p) \cong H_2(\tilde Q_{n,3}; \Z/p)$
consider the equivariant homology spectral sequence for
$G_n$ acting on $\tilde Q_{n,3}$.  

If $p \geq 7$, and a graph $\Gamma_n$ in
$\tilde Q_{n,3}$ has $p$-symmetry, the symmetry
comes from permuting petals attached at the
basepoint and leaving the rest of the graph
fixed.  This is because if an edge $e$ with endpoints
$v$ and $w$, $w\not = *$, is permuted 
nontrivially by the $p$-action, then the $p$-orbit of
$e$ forces the graph to have degree at least 
$p-2$ (if $w$ is fixed by the action) or $p$
(if $w$ is not fixed.)  In any case, the degree would be
too large.  So all $p$-symmetry in $\tilde Q_{n,3}$ 
comes from the rose, in the sense that a graph with $p$-symmetry
consists of a rose graph $R_s$, 
$s \in \{n-6, n-5, \ldots, n\}$,
wedged to some other graph fixed by
$p$.  Such $p$-symmetry results in stabilizers with
comhomology the same as that of symmetric groups.  The homology
of these (see \cite{[N]}) vanishes in dimensions 1 and 2.
Since $H_t(\hbox{stab}_{\Gamma_n}(R_s); \Z/p) = 
H_t(\Sigma_s; \Z/p) = 0$
for $t=1,2$, the result holds by a standard restriction-transfer
argument in group cohomology which implies that we need only be 
concerned with simplices with $p$-symmetry.

If $p=3,5$, there are more simplices with $p$-symmetry in $\tilde Q_{n,3}$,
such as graphs based on $\theta$-graphs
(the most complicated of these in the case $p=3$
being three $\Theta_2$ graphs wedged together
at $*=\circ$ with a symmetry group of 
$(\Sigma_3 \times \Sigma_3 \times \Sigma_3) \rtimes \Sigma_3.$)  
The homology of their stabilizers
still vanishes in dimensions 1 and 2, however, establishing this part
of the lemma.

Mimic the proof of Lemma \ref{lemrat} to obtain the exactness of
$$H_{3}(\tilde Q_{n,3}, \tilde Q_{n,2}; \Z/p) \to
H_2(\tilde Q_{n,2}; \Z/p)
\to H_2(\tilde Q_{n,3}; \Z/p) \to 0.$$ 
\END \end{pf*}

Let $p$ be an odd prime and $n$ be a positive integer.  From the
above lemma, to prove Theorem \ref{lowdim} (1), it suffices to show that
$H_2(\tilde Q_{n,3}; \Z/p) = 0.$  Because part (2) of the theorem will
be established purely by computer program, we illustrate this in some
detail to provide at least one concrete example.

\centerline{\input{holomorph3.pic}}
\begin{figure}[here]
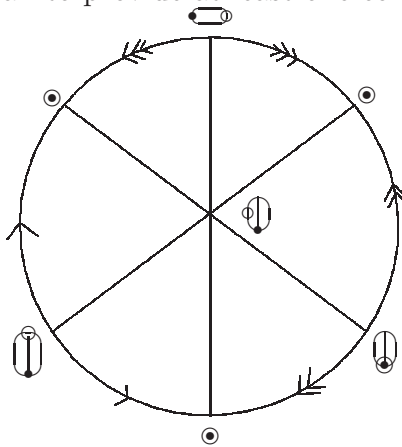

\caption{\label{holomorph3} The $2$-sphere in $\tilde Q_{n,2}$}
\end{figure}

As a notational device, when drawing a graph $\Gamma_{n}$
omit any loops at the basepoint. An example of a $2$-sphere in 
$\tilde Q_{n,2}$ is illustrated
above in Figure \ref{holomorph3}.
In the figure, a filled dot represents the
basepoint $*$ and a hollow dot represents the
other distinguished point $\circ$.  

\centerline{\input{holomorph4.pic}}
\begin{figure}[here]
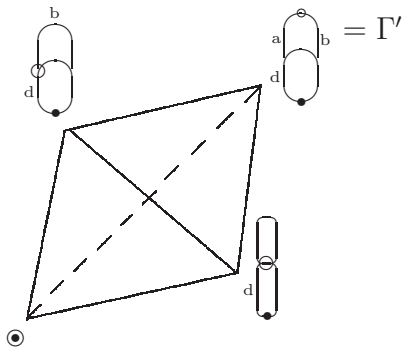

\caption{\label{holomorph4} One of $6$ $3$-simplices making up the sample cube.}
\end{figure}

\begin{figure}[here]
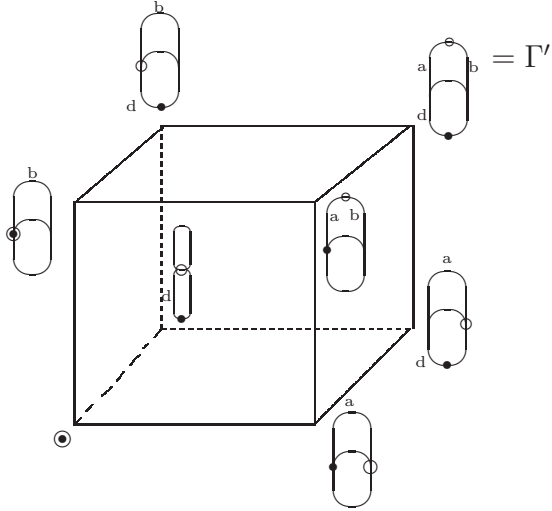

\caption{\label{holomorph5} A sample cube in $\tilde Q_{n,3}$}
\end{figure}

Following \cite{[V]}, we make use of a
cubical structure on some open cells in
$\tilde Q_{n,i}$ and the notion
of plusfaces and minusfaces for this cubical structure and
more generally for the simplicial structure.
Figures \ref{holomorph4} and \ref{holomorph5} provide an example in
the $3$-dimensional complex $\tilde Q_{n,3}$:
Consider a graph $\Gamma'$ with 4 vertices: $*$, $\circ$.
and two other valence $3$ vertices $x$ and $y$.  The graph $\Gamma'$
has $5$ (unoriented) edges.  The edge $a$ connects $\circ$ and $x$, 
$b$ connects $\circ$ and $y$,
$c$ connects $x$ and $y$, $d$ connects $*$ and $x$,
and $e$ connects $*$ and $y$.  The graph $\Gamma'$ is the one
shown in the upper right corners of Figures \ref{holomorph4} and \ref{holomorph5}.
Define a forest $F$ in $\Gamma'$ by $F = \{ a,b,d \}$.  There
are $6$ $3$-simplices in $\tilde Q_{n,2}$ corresponding
to the forest $F$.  Each corresponds to some collapse of
the edges in $F$ in a particular order.  For example, the 
$3$-simplex in Figure \ref{holomorph4} comes from collapsing first $a$,
then $b$, and then $d$.  These $6$ $3$-simplices all fit together
into the cube in Figure \ref{holomorph5}.

\centerline{\input{holomorph5.pic}}

Note that the cube and simplex pictured above each have one
vertex which is maximal (in the poset sense); namely, the one given 
by $\Gamma'$.  A face of the cube or simplex is called a {\em plusface} 
if it is adjacent to this maximal vertex and a {\em minusface} otherwise.
The {\em diagonal} of a cube is the edge (present only in the simplicial
structure and not in the cubical one) joining the poset-maximal and
the poset-minimal vertices of the cube.

The subset $F$ of edges of $\Gamma'$ gives
the interior of a $3$-dimensional cube as
pictured in Figure \ref{holomorph5} because all of the $6$ $3$-simplices 
that form
the cube are distinct.  This in turn is true because no nontrivial
graph automorphism of $\Gamma'$ takes the forest $F$ to itself.
Parts of the boundary of the cube, however, are identified.  For example,
the plusface corresponding to $\{a,b\}$ (the face in back in 
Figure \ref{holomorph5}) ``folds over'' along
its diagonal so that the square forming the plusface is
glued together to form just one $2$-simplex, a triangle.
That is, the square is decomposed into two triangles attached along
a diagonal edge, and in the quotient the two triangles are identified.
The other codimension 1 faces (the two other plusfaces and the three
other minusfaces) of the cube are squares, and
not identified into triangles.

In general, suppose we have a
graph $\Gamma$ which has degree $3$
and a forest
$F = \{ a,b,c \}$ in $\Gamma$.  Then the pair $(\Gamma,F)$ will give
a cube in $\tilde Q_{n}$ if no nontrivial graph automorphism of $\Gamma$ sends
$F$ to itself.  Note that if this is the case, then even though
$(\Gamma,F)$ gives a cube, its faces might be
identified or glued to each other in various ways.
This can happen with both the plusfaces and the
minusfaces.
For example, say $\hat \Gamma$ is the graph obtained from
$\Gamma$ by collapsing the edge $a$.  If a nontrivial
graph automorphism of $\hat \Gamma$ switches $b$ and $c$, the minusface
of the cube corresponding to $\hat \Gamma$ is glued to itself along
a diagonal and is not a square but a $2$-simplex or triangle.

\begin{prop}
\label{claimmod} If $n = 1$ then $\tilde Q_{n,2}$ is contractible.
Otherwise, $\tilde Q_{n,2}$ deformation retracts to the $2$-sphere
pictured in Figure \ref{holomorph3}.
\end{prop}

\begin{pf*}{Proof.}
The case $n=1$ is left to the reader.  Figure \ref{holomorph6} lists
all graphs which give maximal vertices.

\centerline{\input{holomorph6.pic}}
\begin{figure}[here]
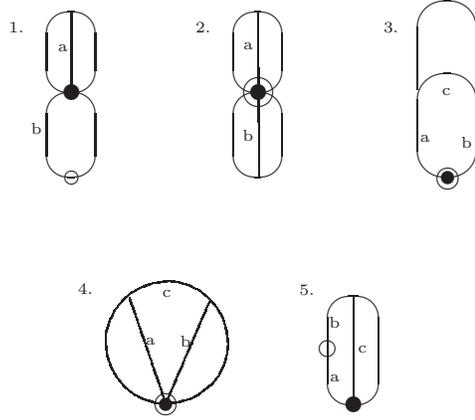

\caption{\label{holomorph6} Graphs giving $2$-simplices}
\end{figure}

Define the deformation retraction onto the $2$-sphere by collapsing
away from the simplices given by the above graphs (numbered 1 - 5) as 
follows:
\begin{enumerate}
\item This graph only has one maximal subforest $\{a,b\}$ so that the
corresponding square (i.e., two $2$-simplices that join together to
form a square) has free plusfaces.  Hence we can collapse this
square away.
\item This graph also only has one maximal subforest $\{a,b\}$
and moreover there is an automorphism of the graph that switches $a$
and $b$.  This graph just contributes one $2$-simplex, the
diagonal of which is automatically a free plusface.
\item The $2$-simplex corresponding to $\{a,b\}$ has a free
diagonal plusface and so can be removed.  In addition, the
square corresponding to $\{a,c\}$ has free plusface $c$.
\item Use exactly the same argument as that for graph $3.$
in Figure \ref{holomorph6}.
\item The three squares that this graph contributes join
together to form the $2$-sphere in Figure \ref{holomorph3}.
\end{enumerate}
\END \end{pf*}

\begin{pf*}{Proof.} {\bf [of Theorem \ref{lowdim} (1).]} From Lemma \ref{lemmod}
and Proposition \ref{claimmod}, we must find an explicit element 
in $H_3(\tilde Q_{n,3}, \tilde Q_{n,2}; \Z/p)$ 
which maps onto the generator of 
$H_2(\tilde Q_{n,2}; \Z/p) = \Z/p.$
The graph $\Gamma'$ from Figure \ref{holomorph5} will be used
to construct the relative cycle.
Basically, the cycle is formed by joining together
the cubes corresponding to the
subforests $\{a,b,d\}$, $\{a,c,d\}$,
and $\{a,c,e\}$.  

The plusface corresponding to $\{a,b\}$ of
$\{a,b,d\}$ folds back onto itself along its
diagonal and so it not free.  The plusface
$\{a,d\}$ of $\{a,b,d\}$ connects up with
the plusface $\{a,d\}$ of $\{a,c,d\}$.  The
remaining plusface $\{b,d\}$ of $\{a,b,d\}$
joins up with the plusface $\{a,e\}$ of
$\{a,c,e\}$.
Moreover, the plusface $\{a,c\}$ of
$\{a,c,d\}$ joins up with the plusface
$\{a,c\}$ of $\{a,c,e\}$.  Finally, the
last remaining plusface $\{c,d\}$ of
$\{a,c,d\}$ connects with the last
remaining plusface $\{c,e\}$ of
$\{a,c,e\}$.

The minusface obtained by collapsing $a$
in $\{a,b,d\}$ is the same as that
obtained by collapsing $a$ in
$\{a,c,d\}$.  The one obtained by
collapsing $b$ in $\{a,b,d\}$ is the
same as what we get if we collapse
$a$ in $\{a,c,e\}$.  The remaining
minusface of $\{a,b,d\}$ from collapsing
$d$ is one of the three squares pictured in
Figure \ref{holomorph3} and corresponds to the
subforest $\{a,b\}$ of graph $5.$ of
Figure \ref{holomorph6}.
Following the same logic, we see that
the minusface of $\{a,c,d\}$ corresponding to
collapsing the edge $c$, is the same
as that of $\{a,c,e\}$ obtained by collapsing
$c$.  The square $\{a,c\}$ of the
graph $5.$ of Figure \ref{holomorph6}
is now seen to be the remaining minusface
acquired from $\{a,c,d\}$ by collapsing $d$.
Similarly, the square $\{b,c\}$ of the graph
$5.$ of Figure \ref{holomorph6} 
is the remaining minusface of $\{a,c,e\}$
which we get if we collapse the edge $e$.

In summary, the cubes $\{a,c,d\}$ 
and $\{a,c,e\}$ glue together to
form a solid 3-ball along the topological 2-disk
formed by the union the plusfaces
$\{c,d\}$ and $\{a,c\}$ of $\{a,c,d\}$ and the
minusface of $\{a,c,d\}$ given by collapsing $c$.
This solid ball is in turn glued to the ball 
corresponding to the cube $\{a,b,d\}$ along the
topological 2-disk formed by the union of the
plusfaces $\{a,d\}$ and $\{b,d\}$ of $\{a,b,d\}$
and the minusfaces of $\{a,b,d\}$ corresponding to
collapsing $\{a\}$ and $\{b\}$.  Note that this 
latter surface is a disk, and not an annulus, because
the plusface $\{a,b\}$ of $\{a,b,d\}$ is
self-identified.  The union of the three cubes is
thus a solid 3-ball with boundary the
2-sphere pictured in Figure \ref{holomorph3}.
\END \end{pf*}

\begin{pf*}{Proof.} {\bf [of Theorem \ref{lowdim} (2).]} 
The methods of \cite{[V]} were used to establish this
by Maple programs.  Briefly, Lemma \ref{lemrat} 
and the cubical structure of $\tilde Q_{n,i}$ are used
to compute the homology groups by enumerating all
relevant graphs and then considering the cubes corresponding
to each graph.
For copies of the specific programs used
and the output files, see: \newline
{\tt http://www.math.uno.edu/\~{}jensen/maple}
\END \end{pf*}

\begin{pf*}{Proof.} {\bf [of Theorem \ref{lowdim} (3).]} 
As mentioned in the introduction, we have a spectral sequence
$$E^2_{r,s} = H_r(Aut(F_n); H_s(F_n; \Q)) \Rightarrow 
H_{r+s}(F_n \rtimes Aut(F_n); \Q).$$ 
The $E^2$-page is 0 except in the row $s=0$, where it is
$H_r(Aut(F_n); \Q)$, and the row $s=1$, where it is
$H_r(Aut(F_n); \Q^n).$  From \cite{[V]}, 
$H_r(Aut(F_n); \Q) = 0$ for $1 \leq r \leq 6$, $n \geq 1$,
except in the case $r=n=4$ where $H_r(Aut(F_n); \Q) = \Q.$ 
Combining this with Theorem \ref{lowdim} (2) yields the desired
result.  The only exceptional cases are where $n=3$ and $n=4$.
When $n=3$, $H_{4}(F_3 \rtimes Aut(F_3); \Q)=\Q$ and
$H_4(Aut(F_3); \Q)=H_5(Aut(F_3); \Q)=0$, forcing  
$H_3(Aut(F_3); \Q^3)=\Q.$  When $n=4$, 
$H_{4}(F_4 \rtimes Aut(F_4); \Q)=H_4(Aut(F_4); \Q)=\Q$ 
and $H_5(Aut(F_4); \Q)=0$, forcing  
$H_3(Aut(F_4); \Q^3)=0.$  
\END \end{pf*}

\end{document}